\RequirePackage{ifpdf}
\ifpdf 
\documentclass[pdftex]{sigma}
\else
\documentclass{sigma}
\fi

\begin{document}

\allowdisplaybreaks

\renewcommand{\thefootnote}{$\star$}

\renewcommand{\PaperNumber}{030}

\FirstPageHeading

\ShortArticleName{Nonlocal Operational Calculi for Dunkl Operators}

\ArticleName{Nonlocal Operational Calculi for Dunkl Operators\footnote{This paper is a contribution to the Special
Issue on Dunkl Operators and Related Topics. The full collection
is available at
\href{http://www.emis.de/journals/SIGMA/Dunkl_operators.html}{http://www.emis.de/journals/SIGMA/Dunkl\_{}operators.html}}}

\newcommand\sign{\mathop{\rm sign}\nolimits}

\Author{Ivan H. DIMOVSKI and Valentin Z. HRISTOV}

\AuthorNameForHeading{I.H. Dimovski and V.Z. Hristov}

\Address{Institute of Mathematics and Informatics, Bulgarian Academy
of Sciences,\\
 Acad. G. Bonchev  Str., Block 8, 1113 Sofia, Bulgaria}

\Email{\href{mailto:dimovski@math.bas.bg}{dimovski@math.bas.bg}, \href{mailto:valhrist@bas.bg}{valhrist@bas.bg}}

\ArticleDates{Received October 15, 2008, in f\/inal form March 04,
2009; Published online March 09, 2009}

\Abstract{The one-dimensional Dunkl operator $D_k$ with a non-negative parameter
$k$, is considered under an arbitrary nonlocal boundary value
condition. The right inverse operator of $D_k$, satisfying this
condition is studied. An operational calculus of Mikusi\'nski type is
developed. In the frames of this operational calculi an extension of
the Heaviside algorithm for solution of nonlocal Cauchy boundary value
problems for  Dunkl functional-dif\/ferential equations  $P(D_k)u=f$
with a given polynomial $P$ is proposed. The solution of these
equations in mean-periodic functions reduces to such problems.
Necessary and suf\/f\/icient condition for existence of unique solution in
mean-periodic functions is found.}


\Classification{44A40; 44A35; 34K06}

\bigskip

\noindent
Here the one-dimensional Dunkl operators $
D_kf(x)={df(x)\over dx}+k{f(x)-f(-x)\over x}$, $k\ge 0$, in
$C^1({\mathbb R})$ under a nonlocal boundary value condition
$\Phi\{f\}=0$ with an arbitrary non-zero linear functional~$\Phi$ in
$C({\mathbb R})$ are considered. The right inverse operators $L_k$ of
$D_k$, def\/ined by $D_kL_kf=f$ and $\Phi\{L_kf\}=0$ are studied. To
this end, the elements of corresponding operational calculi are
developed. A convolution product $f\ast g$ on $C({\mathbb R})$, such
that $L_kf=\{1\}\ast f$, is found. Further, the convolution algebra
$(C({\mathbb R}),\ast)$ is extended to its ring ${\mathfrak M}_k$ of
the multipliers. $(C({\mathbb R}),\ast)$ may be conceived as a part of
${\mathfrak M}_k$ due to the embedding $f\hookrightarrow f\ast$. The
ring ${\cal M}_k$ of multiplier fractions~$ {A\over B}$,
such that $A,B\in{\mathfrak M}_k$ and $B$ being non-divisor of zero in
the operator multiplication, is constructed.

A Heaviside algorithm for ef\/fective solution of nonlocal Cauchy
boundary value problems for {\it Dunkl functional-differential
equations} $P(D_k)u=f$ with polynomials $P$ is developed. The
solution of these equations in mean-periodic functions reduces to such
problems. Necessary and suf\/f\/icient condition for existence of unique
solution in mean-periodic functions is found.

The operational calculus, developed here, is a generalization of the
nonlocal operational calculus for $ D_0={d\over dx}$ (see
Dimovski~\cite{Di2}). Some background material about the Dunkl
operators is taken from our previous paper~\cite{D-H-S} without proofs.

\section[The right inverse operators of $D_k$ in $C({\mathbb R})$ and
corresponding Taylor formulae]{The right inverse operators of $\boldsymbol{D_k}$ in $\boldsymbol{C({\mathbb R})}$\\ and corresponding Taylor formulae}\label{2}

Let $L_k$ denote an arbitrary right inverse operator of $D_k$ in
$C({\mathbb R})$. First, we consider a special right inverse
$\Lambda_k$ of $D_k$, where $y(x)=\Lambda_kf(x)$ for $f\in C({\mathbb
R})$ is the solution of the equation $D_ky=f(x)$ with initial
condition $y(0)=0$.

\begin{lemma}\label{l1}
The right inverse operator $\Lambda_k$ of $D_k$, defined by the
initial condition $\Lambda_kf(0)=0$ has the form
\begin{equation*}
\Lambda_kf(x)=\int_0^x\left[f_{\rm o}(t)+\left({t\over
x}\right)^{2k}f_{\rm e}(t)\right]dt,
\end{equation*}
where $f_{\rm e}$ and $f_{\rm o}$ are the even and the odd parts of $f$,
respectively.
\end{lemma}
The proof is a matter of a simple check (see~\cite[p.~198]{D-H-S}).

In the general case, an arbitrary right inverse operator $L_k$ of $D_k$
has a representation of the form
\[
L_kf(x)=\int_0^x\left[f_{\rm e}(t)+\left({t\over
x}\right)^{2k}f_{\rm o}(t)\right]dt+C.
\]
In order $L_k$ to be a linear operator, the additive constant $C$
should
depend on $f$ and to be a~linear functional $\Psi\{f\}$ in~$C({\mathbb
R})$. Hence, an arbitrary linear right inverse operator $L_k$ of~$D_k$
in~$C({\mathbb R})$ has the form
\[
L_kf(x)=\Lambda_kf(x)+\Psi\{f\},
\]
with a linear functional $\Psi$ in $C({\mathbb R})$.

According to the general theory of right invertible operators (Bittner~\cite{Bi}, Przeworska-Ro\-le\-wicz~\cite{Pr}), an important characteristic
of $L_k$ is its {\it initial projector}
\begin{equation}\label{inproj}
Ff(x)=f(x)-L_kD_kf(x)=\Phi\{f\}.
\end{equation}
It maps $C^1({\mathbb R})$ onto $\ker D_k={\mathbb C}$, i.e.\ it is a
linear functional $\Phi$ on $C^1({\mathbb R})$. This identity written
in the form
\begin{equation}\label{inproj1}
L_kD_kf(x)=f(x)-\Phi\{f\}.
\end{equation}
will be used later.
Expressing $\Phi$ by $\Psi$, we obtain
\[
\Phi\{f\}=f(0)-\Psi\{D_kf\}.
\]
Let us note that $\Phi\{1\}=1$, which expresses the projector property
of~$F$.

Considering the right inverse operator $L_k$ of $D_k$, it is more
convenient to look on $L_kf=y$ as the solution of an elementary
boundary value problem of the form
\begin{equation*}
D_ky=f,\qquad\Phi\{y\}=0,
\end{equation*}
assuming that $\Phi$ is a given linear functional on $C({\mathbb R})$
with $\Phi\{1\}=1$. This restriction of the class of right inverse
operators $L_k$ of $D_k$ is adequate when we are to consider nonlocal
Cauchy problems for Dunkl equations.

\begin{theorem}
Let $\Phi:C({\mathbb R})\to{\mathbb C}$ be a linear functional, such
that $\Phi\{1\}=1$. Then the right inverse operator $L_k$ of $D_k$,
defined by the boundary value condition $\Phi\{L_kf\}=0$ has the form
\begin{equation*}
L_kf(x)=\int_0^x\left[f_{\rm e}(y)+\left({y\over
x}\right)^{2k}f_{\rm o}(y)\right]dy-
\Phi_t\left\{\int_0^t\left[f_{\rm e}(y)+\left({y\over
t}\right)^{2k}f_{\rm o}(y)\right]dy\right\}.
\end{equation*}
\end{theorem}
The proof follows immediately from Lemma \ref{l1} and the condition
$\Phi\{1\}=1$.
\begin{definition}\label{d1}
The polynomials
\begin{equation}\label{e1}
A_{k,n}(x)=L_k^n\{1\}(x),\qquad n=0,1,2,\dots
\end{equation}
are said to be {\it Dunkl--Appell polynomials}.
\end{definition}

\begin{lemma}\label{l2}
The {\it Dunkl--Appell polynomials system}
$\{A_{k,n}(x)\}_{n=0}^{\infty}$ satisfies the recurrences
\begin{equation}\label{e2}
A_{k,0}(x)\equiv 1,\qquad{\rm and}\qquad
D_kA_{k,n+1}(x)=A_{k,n}(x),\qquad\Phi\{A_{k,n+1}\}=0,\qquad
n\ge 0
\end{equation}
and conversely, \eqref{e2} implies \eqref{e1}.
\end{lemma}
The check is immediate. Similar polynomials are introduced implicitly
by M.~R\"osler and M.~Voit \cite[p.~346]{R-V}.

\begin{lemma}[Taylor formula with remainder term] If $f\in
C^{(N)}({\mathbb R})$, then
\begin{equation}\label{macl}
f(x)=\sum_{j=0}^{N-1}\Phi\big\{D_k^jf\big\}A_{k,j}(x)+L_k^N\big(D_k^Nf\big)(x),
\end{equation}
where $A_{k,j}(x)=L_k^j\{1\}(x)$ are Dunkl--Appell polynomials.
\end{lemma}
This formula is an analogue of the particular case of the Taylor
formula known as the Maclaurin formula.

\begin{proof}
Delsarte~\cite{De}, Bittner~\cite{Bi}, and~Przeworska-Rolewicz
\cite{Pr} give variants of the Taylor formula for right invertible
operators in linear spaces. In our case (\ref{macl}) can be written as
\[
I=\sum_{j=0}^{N-1}L_k^jFD_k^j+L_k^ND_k^N,
\]
where $I$ is the identity operator and $F=I-L_kD_k$.
In functional form the above identity takes the form
\[
f(x)=\sum_{j=0}^{N-1}L_k^jFD_k^jf(x)+L_k^ND_k^Nf(x),
\]
where the initial projector $F$ of $L_k$ (\ref{inproj}) is the linear
functional $\Phi$:
\[
Ff(x)=f(x)-L_kD_kf(x)=\Phi\{f\}.
\]
$F$ projects the space $C({\mathbb R})$ onto the space ${\mathbb C}$ of
the constants. Hence
\begin{equation*}
f(x)=\sum_{j=0}^{N-1}\Phi\big\{D_k^jf\big\}L_k^j\{1\}(x)+L_k^ND_k^Nf(x),
\end{equation*}
which is the Taylor formula (\ref{macl}).
\end{proof}

\section[Convolutional products for the right inverses $L_k$ of
$D_k$]{Convolutional products for the right inverses $\boldsymbol{L_k}$ of
$\boldsymbol{D_k}$}\label{3}

In Dunkl \cite[Theorem 5.1]{Du} the similarity
operator
\begin{equation*}
V_kf(x)=b_k\int_{-1}^1f(xy)(1-y)^{k-1}(1+y)^kdy,\qquad
b_k={\Gamma(2k+1)\over 2^{2k}\Gamma(k)\Gamma(k+1)}
\end{equation*}
is found, which transforms the dif\/ferentiation operator
$ D={d\over dx}$ into $D_k$:
\[
V_kD=D_kV_k.
\]
Usually this operator is called {\it intertwining operator}. The
constant $b_k$ is chosen to ensure that $V_k\{1\}=1$.

The problem of inverting the Dunkl intertwining operator $V_k$ is
discussed by several authors, see e.g.\ Trim\`eche~\cite{T}, Betankor,
Sif\/i, Trim\`eche~\cite{B-S-T}, but we will use the explicit formulae
from Ben~Salem and Kallel~\cite[p.~159]{S-K}.

Denoting $ Sf(x)={1\over 2x}{df(x)\over dx}$, the inverse
$V_k^{-1}$ of $V_k$ has the following representations:

(i) If $k=n+r$ is non-integer with integer part $n$ and
$r\in(0,1)$, then
\begin{gather*}
V_k^{-1}f(x) =
c_k\left[
|x|S^{n+1}\left\{\int_0^{|x|}\big(x^2-y^2\big)^{-r}f_{\rm e}(y)y^{2k}dy\right\}
\right.\\
\left.\phantom{V_k^{-1}f(x) =}{} +
\sign(x)S^{n+1}\left\{\int_0^{|x|}\big(x^2-y^2\big)^{-r}f_{\rm o}(y)y^{2k+1}
dy\right\}
\right],\qquad x\ne 0,
\end{gather*}
where
$  c_k={2\sqrt{\pi}\over\Gamma\left(n+r+{1\over 2}\right)
\Gamma(1-r)}$.

(ii) If $k$ is a non-negative integer, then
\[
V_k^{-1}f(x)={\sqrt{\pi}\over\Gamma\left(k+{1\over 2}\right)}
\big[xS^k\big(x^{2k-1}f_{\rm e}(x)\big)+S^k\big(x^{2k}f_{\rm o}(x)\big)\big],\qquad x\ne 0.
\]

$V_k$ transforms $C({\mathbb R})$ into a proper subspace
$\widetilde{C_k}=V_k(C({\mathbb R}))$ of it. $V_k$ is a similarity from
a~right inverse operator $\Lambda$ of $\displaystyle D_0={d\over dx}$
to $L_k$. In order to specify the operator $\Lambda$ let us def\/ine the
linear functional
\begin{equation*}
\widetilde{\Phi}\{f\}=(\Phi\circ V_k)\{f\}
\end{equation*}
in $\widetilde{C_k}$. Then def\/ine
$\Lambda:\widetilde{C_k}\to\widetilde{C_k}$ to be the solution
$y=\Lambda\widetilde{f}$ of the elementary boundary value problem
\[
D_0y(x)\equiv y'(x)=\widetilde{f}(x),\qquad\widetilde{\Phi}\{y\}=0.
\]
This solution has the form
\begin{equation*}
\Lambda\widetilde{f}(x)=\int_0^x\widetilde{f}(y)dy-
\widetilde{\Phi_{t}}\left\{\int_0^t\widetilde{f}(\tau)d\tau\right\}.
\end{equation*}

\begin{lemma}\label{simil}
The following similarity relation holds
\[
V_k\Lambda=L_kV_k.
\]
\end{lemma}

\begin{proof}
Applying $V_k$ to the def\/ining equation
$D(\Lambda\widetilde{f})=\widetilde{f}$, one obtains
\[
V_kD(\Lambda\widetilde{f})=V_k\widetilde{f}=f
\qquad \mbox{or}\qquad
D_k(V_k\Lambda\widetilde{f})=V_k\widetilde{f}=f.
\]
In fact, the boundary value condition
$\widetilde{\Phi}\{\Lambda\widetilde{f}\}=0$ can be written as
$\Phi\{V_k\Lambda\widetilde{f}\}=0$. Hence $u=V_k\Lambda\widetilde{f}$
is the solution of the boundary value problem $D_ku=f$, $\Phi\{u\}=0$,
i.e.~$u=L_kf$. Therefore
\[
V_k\Lambda V_k^{-1}f=L_kf\qquad\text{or}\qquad V_k\Lambda=L_kV_k.
\tag*{\qed}
\]
\renewcommand{\qed}{}
\end{proof}

The similarity relation (\ref{simil}) allows to
introduce a convolution structure $\ast:C({\mathbb R})\times
C({\mathbb R})\to C({\mathbb R})$, such that $L_k$ to be the
convolution operator $L_k=\{1\}\ast$ in $C({\mathbb R})$.

The operator $\Lambda$ is def\/ined not only in $\widetilde{C_k}$, but in
the whole space $C({\mathbb R})$. This allows to introduce a
convolution structure $\widetilde{\ast}:C({\mathbb R})\times
C({\mathbb R})\to C({\mathbb R})$.

\begin{lemma}\label{convt}
The operation
\begin{equation}\label{convtilde}
(\widetilde{f} \,\widetilde{\ast}\, \widetilde{g})(x)=
\widetilde{\Phi}_t
\left\{\int_t^x\widetilde{f}(x+t-\tau)\widetilde{g}(\tau)d\tau\right\}
\end{equation}
is a bilinear, commutative and associative operation in
$\widetilde{C_k}=V_k(C({\mathbb R}))$ such that
\begin{equation}\label{convrepr}
\Lambda\widetilde{f}=\{1\}\, \widetilde{\ast}\, \widetilde{f}.
\end{equation}
It satisfies the boundary value condition
$\widetilde{\Phi}\{\widetilde{f}\,\widetilde{\ast}\, \widetilde{g}\}=0$ for
arbitrary $\widetilde{f}$ and $\widetilde{g}$ in $C({\mathbb R})$.
\end{lemma}

The proof of the assertion that
$\widetilde{f}\, \widetilde{\ast}\, \widetilde{g}$ is an inner operation in
$\widetilde{C_k}$ follows directly from the explicit
inversion formula for $V_k$ (see Xu~\cite{Xu} or Ben~Salem and Kallel
\cite[Theorem~1.1]{S-K}). In Dimovski \mbox{\cite[p.~52]{Di}} it is proved
that (\ref{convtilde}) is a bilinear, commutative and associative
operation in~$C({\mathbb R})$, and hence in
$\widetilde{C_k}=V_k(C({\mathbb R}))$. The second relation~(\ref{convrepr}) is obvious. The proof of
$\widetilde{\Phi}\{\widetilde{f}\widetilde{\ast}\widetilde{g}\}=0$ is
also elementary (see Dimovski \cite[p.~54]{Di}).

\begin{theorem}\label{conv}
The operation
\begin{equation}\label{convol}
f\ast g=D_k^{2n}V_k\big[\big(V_k^{-1}L_k^nf\big) \, \widetilde{\ast}\,
\big(V_k^{-1}L_k^ng\big)\big],
\end{equation}
where $n$ is the integer part of $k$, is a convolution of $L_k$ in
$C({\mathbb R})$ such that
\begin{equation}\label{conv2}
L_kf=\{1\}\ast f
\end{equation}
and the boundary value condition $\Phi\{f\ast g\}=0$ is satisfied for
arbitrary $f$ and $g$ in $C({\mathbb R})$.
\end{theorem}

\begin{proof}
The assertion of the theorem follows from Lemmas \ref{convt} and
\ref{simil} and a general theorem of Dimovski~\cite[Theorem~1.3.6, p.~26]{Di}. This convolution is introduced in Dimovski, Hristov and Sif\/i~\cite{D-H-S}.
\end{proof}

\begin{remark}
The convolution (\ref{convol}) reduces to
\[
f\ast g=V_k\big[\big(V_k^{-1}f\big)\,\widetilde{\ast}\,\big (V_k^{-1}g\big)\big]
\]
for $n=0$, i.e.\ when $0<k<1$.
\end{remark}

From (\ref{conv2}) and Def\/inition~\ref{d1} it follows that
\[
L_k^{N+1}f=\{A_{k,N}\}\ast f,
\]
where $A_{k,N}$ is the Dunkl--Appell polynomial of degree exactly $N$.
This allows also to state the Taylor formula~(\ref{macl}) with
remainder term in the Cauchy form:

\begin{lemma}
If $f\in C^{(N)}({\mathbb R})$, then
\[
f(x)=\sum_{j=0}^{N-1}\Phi\big\{D_k^jf\big\}A_{k,j}(x)+\big(A_{k,N-1}\ast
D_k^Nf\big)(x),
\]
where $A_{k,j}(x)$, $j=0,1,2,\dots,N-1$, are the Dunkl--Appell
polynomials $A_{k,j}(x)=L_k^j\{1\}$.
\end{lemma}

\section[The ring of multipliers of the convolutional algebra
$(C({\mathbb R}),\ast)$]{The ring of multipliers of the convolutional algebra
$\boldsymbol{(C({\mathbb R}),\ast)}$}\label{5}

The convolutional algebras $(C({\mathbb R}),\ast)$ with convolution
product (\ref{convol}), are annihilators-free (or algebras without
order in the terminology of Larsen \cite[p.~13]{L}). This means that
in each of these algebras $f\ast g=0$, $\forall \, g\in C({\mathbb R})$,
implies $f=0$.

\begin{definition}\label{multipl}
An operator $A:C({\mathbb R})\to C({\mathbb R})$ is said to be a {\it
multiplier of the convolutional algebra} $(C({\mathbb R}),\ast)$ if\/f
\begin{equation}\label{defmult}
A(f\ast g)=(Af)\ast g
\end{equation}
for arbitrary $f,g\in C({\mathbb R})$.
\end{definition}

As it is shown in Larsen \cite{L}, it is not necessary to assume
neither that $A$ is a linear operator, nor that it is continuous in
$C({\mathbb R})$. These properties of the multipliers follow
automatically from~(\ref{defmult}). Something more, a general result
of Larsen \cite[p.~13]{L} implies

\begin{theorem}\label{commring}
The set of the multipliers of the convolutional algebra $(C({\mathbb
R}),\ast)$ form a commutative ring $\mathfrak{M}_k$.
\end{theorem}

The simplest multipliers of $(C({\mathbb R}),\ast)$ are the numerical
operators $[\alpha]$ for $\alpha\in{\mathbb C}$, def\/ined by
\begin{equation*}
[\alpha]f=\alpha f,\qquad\forall \,f\in C({\mathbb R}),
\end{equation*}
and the convolutional operators $f\ast\ $ for $f\in C({\mathbb R})$,
def\/ined by
\begin{equation*}
(f\ast)g=f\ast g,\qquad\forall \,g\in C({\mathbb R}).
\end{equation*}

Further we need the following characterization result for the
multipliers of $(C({\mathbb R}),\ast)$:

\begin{theorem}\label{charact}
A linear operator $A:C({\mathbb R})\to C({\mathbb R})$ is a multiplier
of $(C({\mathbb R}),\ast)$ iff it admits a~representation of the form
\begin{equation}\label{AD}
Af=D_k(m\ast f),
\end{equation}
where the function $m=A\{1\}$ is such that $m\ast f\in C^1({\mathbb
R})$ for all $f\in C({\mathbb R})$.
\end{theorem}

\begin{proof}
Let $A:C({\mathbb R})\to C({\mathbb R})$ be a multiplier of
$(C({\mathbb R}),\ast)$. The operator $L_kf=\{1\}\ast f$ is also a
multiplier. Then, according to Theorem~\ref{commring},
\[
AL_k=L_kA.
\]
Applying $A$ to $L_kf=\{1\}\ast f$, we get
\[
L_kAf=AL_kf=A(\{1\}\ast f)=(A\{1\})\ast f.
\]
The identity
\begin{equation}\label{LA}
L_k(Af)=m\ast f
\end{equation}
with $m=A\{1\}$ is possible only if $m\ast f\in C^1({\mathbb R})$ for
each $f\in C({\mathbb R})$. It remains to apply $D_k$ to~(\ref{LA}) in
order to obtain~(\ref{AD}).

Conversely, let $A:C({\mathbb R})\to C({\mathbb R})$ be the operator
def\/ined by (\ref{AD}), i.e.\ $Af=D_k(m\ast f)$, where $m\in C({\mathbb
R})$ is such that $m\ast f\in C^1({\mathbb R})$ for all $f\in
C({\mathbb R})$. Then
\[
A(f\ast g)=D_k(m\ast(f\ast g))=D_k((m\ast f)\ast g).
\]
But $m\ast f=L_kD_k(m\ast f)$ due to formula (\ref{inproj1}) since
$\Phi(m\ast f)=0$ by Theorem~\ref{conv}. Then
\[
A(f\ast g)=D_kL_k[D_k(m\ast f)\ast g]=(Af)\ast g.
\]
Hence $A$ is a multiplier of the convolution algebra $(C({\mathbb
R}),\ast)$.
\end{proof}

The specif\/ication of the function $m=A\{1\}$ is, in general, a
nontrivial problem even in the case of the simplest Dunkl operator
$  D_0={d\over dx}$ (the usual dif\/ferentiation). This
could be conf\/irmed by the following two examples:

\begin{example}
If $\Phi\{f\}=f(0)$, then $m$ is a continuous function of locally
bounded variation, i.e.\ $m\in BV\cap C({\mathbb R})$ (see Dimovski
\cite[p.~26]{Di}).
\end{example}

\begin{example}
Let $\Phi\{f\}=\int_0^1f(x)dx$. Then $m\in C({\mathbb R})$ can be
arbitrary (see Dimovski \cite[p.~69]{Di}).
\end{example}

\section[Nonlocal operational calculi for $D_k$]{Nonlocal operational calculi for $\boldsymbol{D_k}$}\label{6}

Our aim here is to develop a direct operational calculus for solution
of the following {\it nonlocal Cauchy problem} for the operator~$D_k$:
Solve the equation $P(D_k)u=f$ with a polynomial $P$ and a given $f\in
C({\mathbb R})$ under the boundary value conditions
$\Phi\{D_k^ju\}=\alpha_j$, $j=0,1,2,\dots,\deg P-1$, where $\alpha_j$
are given constants and $\Phi$ is a nonzero linear functional on
$C({\mathbb R})$.

This is a special case of the problems considered by R.~Bittner
\cite{Bi} and D.~Przeworska-Rolewicz~\cite{Pr} for an arbitrary right
invertible operator $D$ instead of~$D_k$.

Our intention here is to propose constructive results and to obtain an
explicit solution of the boundary value problems considered. This is
done by means of an operational calculus essential part of which is an
extension of the Heaviside algorithm.

This operational calculus is developed using a direct algebraic
approach based on the convolution (\ref{convol}). Instead of
Mikusi\'nski's method~\cite{M} of convolutional fractions
$ {f\over g}$, we follow an alternative approach of
multiplier fractions $ {A\over B}$, where $A$ and $B$ are
multipliers of the convolutional algebra $(C({\mathbb R}),\ast)$ and $B$
is a non-divisor of zero in the operator multiplication.

Let us consider the ring $\mathfrak{M}_k$ of the multipliers of the
convolutional algebra $(C({\mathbb R}),\ast)$. The correspondence
$\alpha\mapsto[\alpha]$ is an embedding of ${\mathbb C}$ into
$\mathfrak{M}_k$. The correspondence $f\mapsto f\ast$ is an embedding
of $(C({\mathbb R}),\ast)$ in $\mathfrak{M}_k$. Hence, we may
consider ${\mathbb C}$ and $C({\mathbb R})$ as parts of
$\mathfrak{M}_k$.

$\mathfrak{M}_k$ is a commutative ring (Theorem~\ref{commring}). The
subset $\mathfrak{N}_k$ of $\mathfrak{M}_k$, consisting of the
non-zero non-divisors of zero with respect to the operator
multiplication in $\mathfrak{M}_k$, is nonempty. Indeed, at least the
identity operator $I$ and the right inverse $L_k$ of $D_k$ belong to
$\mathfrak{N}_k$. In addition, $\mathfrak{N}_k$ is a~multiplicative
subset, i.e.\ if $A,B\in\mathfrak{N}_k$, then $AB\in\mathfrak{N}_k$.

Consider the Cartesian product
\[
\mathfrak{M}_k\times\mathfrak{N}_k=
\{(A,B):A\in\mathfrak{M}_k,B\in\mathfrak{N}_k\}
\]
and introduce the equivalence relation
\begin{equation}\label{equivrel}
(A,B)\sim(A',B') \ \Leftrightarrow \ AB'=BA'.
\end{equation}

\begin{definition}
The set ${\cal M}_k=
\mathfrak{M}_k\times\mathfrak{N}_k/_{\displaystyle\sim}$
obtained by the factorization of $\mathfrak{M}_k\times\mathfrak{N}_k$
with respect to the equivalence relation {\rm(\ref{equivrel})} is said
to be the {\it ring of multiplier fractions}.
\end{definition}

${\cal M}_k$ may be considered both as an extension of the f\/ield
${\mathbb C}$ of the complex numbers and of the ring $(C({\mathbb
R}),\ast)$. Formally, this is seen by the embeddings
\[
\alpha\hookrightarrow{[\alpha]\over I} \qquad\textrm{and}\qquad
f\hookrightarrow{f\ast\over I}.
\]
In the sequel we denote the identity operator $I$ simply by $1$. The
multiplication operation of the two elements $p$ and $q$ in ${\cal
M}_k$ will be denoted simply by $pq$. Therefore, instead of $f\ast g$
we will write $fg$.

For our aims the most important elements of ${\cal M}_k$ are
\[
L_k=\{1\}\qquad\textrm{and}\qquad S_k={1\over L_k}.
\]
The fraction $S_k$ with the identity operator as numerator and with
$L_k$ as denominator will be called {\it algebraic Dunkl operator}.
Its relation to the ordinary Dunkl operator $D_k$ is given by the
following theorem:

\begin{theorem}\label{SkDk}
Let $f\in C^1({\mathbb R})$. Then
\begin{equation}\label{DSPhi}
D_kf=S_kf-\Phi\{f\}.
\end{equation}
\end{theorem}
Note that identity {\rm(\ref{DSPhi})} should be interpreted as
\[
(D_kf)\ast=S_k(f\ast)-[\Phi\{f\}],
\]
where $(D_kf)\ast$ and $(f\ast)$ are to be understood as convolution
operators and $[\Phi\{f\}]$ as the numerical operator determined by
the number $\Phi\{f\}$. $S_k$ is neither convolutional nor numerical
operator, but an element of ${\cal M}_k$.
\begin{proof}
In Section \ref{2} (equality (\ref{inproj1})) we have seen that
\[
L_kD_kf=f-\Phi\{f\},
\]
where $\Phi\{f\}$ is the corresponding constant function
$\{\Phi\{f\}\}$. Considered as an operator identity, this can be
written as $(L_kD_kf)\ast=f\ast-\{\Phi\{f\}\}\ast$ or
$L_k[D_k(f\ast)]=f\ast-\Phi\{f\}  L_k$. Hence
\[
L_k(D_kf)\ast=(f\ast)-\Phi\{f\}L_k.
\]
It remains to multiply by $S_k$ to obtain (\ref{DSPhi}).
\end{proof}

Relation (\ref{DSPhi}) may be characterized as the basic formula of
our operational calculus. Using it repeatedly, we obtain

\begin{corollary}
Let $f\in C^{(N)}({\mathbb R})$. Then
\begin{equation}\label{Dkn}
D_k^Nf=S_k^Nf-\sum_{j=0}^{N-1}\Phi\big\{D_k^jf\big\}S_k^{N-j-1}.
\end{equation}
\end{corollary}

\begin{remark}
The last formula is equivalent to the Taylor formula (\ref{macl}) in
Section \ref{2}.
\end{remark}
Let
$P(\lambda)=a_0\lambda^m+a_1\lambda^{m-1}+\cdots+a_{m-1}\lambda+a_m$,
$a_0\ne 0$, and $\Phi$ be a non-zero linear functional on $C({\mathbb
R})$.
\begin{definition}\label{nlcp}
The problem for solving the Dunkl functional-dif\/ferential equation
\[
P(D_k)u=f,\qquad f\in C({\mathbb R})
\]
under the boundary value conditions
\[
\Phi\big\{D_k^ju\big\}=\alpha_j,\qquad j=0,1,2,\dots,m-1
\]
is called a {\it nonlocal Cauchy problem determined by the functional}
$\Phi$.
\end{definition}

By means of (\ref{DSPhi}) and (\ref{Dkn}) it is possible to
``algebraize'' any nonlocal Cauchy boundary value problem.

The simplest nonlocal Cauchy problem for $D_k$,
determined by a linear functional $\Phi$ in $C({\mathbb R})$ concerns
the functional-dif\/ferential equation
\begin{equation*}
D_ku(x)-\lambda u(x)=f(x)
\end{equation*}
with the boundary condition $\Phi\{u\}=0$.

It is known that the solution of the homogeneous equation
\[
D_ku(x)-\lambda u(x)=0
\]
under the initial condition $u(0)=1$ is
\begin{equation*}
u_k(\lambda x)=j_{k-{1\over 2}}(i\lambda x)+{\lambda x\over
2k+1}j_{k+{1\over 2}}(i\lambda x)
\end{equation*}
(see Ben~Salem and Kallel \cite[p.~161]{S-K}), where $j_{\alpha}(x)$
denotes the modif\/ied (normalized) Bessel function
\[
j_{\alpha}(x)=2^{\alpha}\Gamma(\alpha+1){J_{\alpha}(x)\over
x^{\alpha}} ,\qquad x\ne 0 \qquad \textrm{and} \qquad j_{\alpha}(0)=1.
\]

We introduce the {\it Dunkl indicatrix} of the functional $\Phi$ as
the following entire function of exponential type:
\begin{equation*}
E_k(\lambda)=\Phi_{\xi}\{u_k(\lambda\xi)\}=
\Phi_{\xi}\left\{j_{k-{1\over 2}}(i\lambda\xi)+{\lambda\xi\over
2k+1}j_{k+{1\over 2}}(i\lambda\xi)\right\}.
\end{equation*}

\begin{lemma}
The function ${u_k(\lambda x)\over E_k(\lambda)}$ is the
generating function of the Dunkl--Appell polynomials system, i.e.
\[
{u_k(\lambda x)\over
E_k(\lambda)}=\sum_{n=0}^{\infty}\lambda^nA_{k,n}(x).
\]
\end{lemma}

Here we will skip the simple proof.
The linear operator $L_{k,\lambda}$ def\/ined as the solution
$u(x)=L_{k,\lambda}f(x)$ of the nonlocal Cauchy boundary value problem
\[
D_ku-\lambda u=f ,\qquad \Phi\{u\}=0,
\]
is said to be the {\it resolvent operator of the Dunkl operator} under
the boundary value condition $\Phi\{u\}=0$.

\begin{theorem}\label{resolvent}
The resolvent operator $L_{k,\lambda}$ admits the convolutional
representation
\begin{equation*}
L_{k,\lambda}f(x)=l_k(\lambda,x)\ast f(x),\qquad{\rm where}\qquad
l_k(\lambda,x)={u_k(\lambda x)\over E_k(\lambda)}.
\end{equation*}
\end{theorem}

\begin{proof}
We will use the formula
\[
D_k(f\ast g)=(D_kf)\ast g+\Phi\{f\}g
\]
which is true under the assumption $f\in C^1({\mathbb R})$. It follows
from a more general result of Dimovski \cite[Theorem~1.38]{Di}, but in our
case it can be verif\/ied directly. It gives
\[
D_k\{l_k(\lambda,x)\ast f(x)\}=
D_kl_k(\lambda,x)\ast f(x)+
\Phi_{\xi}\{l_k(\lambda,\xi)\}f(x)=
\lambda\{l_k(\lambda,x)\ast f(x)\}+f(x).
\]
Hence $\displaystyle u=\{l_k(\lambda,x)\ast f(x)\}$ satisf\/ies the
equation $D_ku-\lambda u=f$. It remains to verify the boundary value
condition $\Phi\{u\}=0$. But it follows from the basic property
$\Phi\{f\ast g\}=0$ of the convolution (Theorem \ref{conv}).
\end{proof}

The resolvent operator $L_{k,\lambda}$ exists for each $\lambda$ with
$E_k(\lambda)\ne 0$. The zeros of $E_k(\lambda)$ are the eigenvalues
of the boundary value problem $D_ku-\lambda u=0$, $\Phi\{u\}=0$. They
form  an enumerable set
$\{\lambda_1,\lambda_2,\dots,\lambda_n,\dots\}$ except in the case
when $\Phi$ is a Dirac functional $\Phi\{f\}=f(a)$, when
$E_k(\lambda)\ne 0$ for all $\lambda\in{\mathbb C}$.

It is easy to f\/ind the solution of our problem in ${\cal M}_k$. Using
the basic formula of the operational calculus (see Theorem
\ref{SkDk}), we have $D_ku=S_ku$ since $\Phi\{u\}=0$, and then
\[
S_ku-\lambda u=f\qquad{\rm or}\qquad (S_k-\lambda)u=f.
\]
In order to write the solution
\[
u={1\over S_k-\lambda}f
\]
we must be sure that $S_k-\lambda$ is non-divisor of zero.

\begin{lemma}\label{nondiv}
$S_k-\lambda$ is a divisor of zero in ${\cal M}_k$ iff
$E_k(\lambda)=0$.
\end{lemma}
\begin{proof}
Let $S_k-\lambda$ be a divisor of zero in ${\cal M}_k$. Then there
exists a multiplier fraction ${A\over B}$ such that $A\ne
0$ and
\[
(S_k-\lambda){A\over B}=0,
\]
which is equivalent to $(S_k-\lambda)A=0$. Since $A\ne 0$, then there
is a function $g\in C({\mathbb R})$ such that $Ag=v\ne 0$. Then
\[
(S_k-\lambda)v=0.
\]
Multiplying by $L_k$ we get
\[
(1-\lambda L_k)v=0\qquad{\rm or}\qquad v-\lambda L_kv=0.
\]
Since $\Phi(L_kv)=0$ by the def\/inition of $L_k$ (Section~\ref{2}),
then $\Phi\{v\}=0$.

Applying $D_k$, we get $D_kv-\lambda v=0$, $\Phi\{v\}=0$. According to
Ben~Salem and Kallel \cite{S-K}, all the non-zero solutions of
$D_kv-\lambda v=0$ are $v=C(j_{k-{1\over 2}}(i\lambda x)+{\lambda
x\over 2k+1}j_{k+{1\over 2}}(i\lambda x))$ with a constant $C\ne 0$.
The boundary value condition $\Phi\{v\}=0$ is equivalent to
$E_k(\lambda)=0$.

Conversely, if $E_k(\lambda)=0$, then there exists a solution $v\ne 0$
of the eigenvalue problem $D_kv-\lambda v=0$, $\Phi\{v\}=0$. For this
$v$ we have
\[
(S_k-\lambda)v=0
\]
and hence $S_k-\lambda$ is a divisor of zero in ${\cal M}_k$.
\end{proof}

\begin{theorem}\label{1sk}
Let $\lambda\in{\mathbb C}$ be such that $E_k(\lambda)\ne 0$. Then
\begin{equation}\label{skl}
{1\over S_k-\lambda}=\{l_k(\lambda,x)\}\ast=
{1\over E_k(\lambda)}\left\{j_{k-{1\over 2}}(i\lambda x)+{\lambda
x\over 2k+1}j_{k+{1\over 2}}(i\lambda x)\right\}\ast.
\end{equation}
\end{theorem}
\begin{proof}
We have seen that
\[
L_{k,\lambda}f(x)=\{l_k(\lambda,x)\}\ast f.
\]
But for the solution $u=L_{k,\lambda}f$ of the boundary value problem
$D_ku-\lambda u=f$, $\Phi\{u\}=0$, in the case $E_k(\lambda)\ne 0$ we
found
\[
u={1\over S_k-\lambda}f.
\]
Since the convolution $\ast$ is annihilators-free, then (\ref{skl})
follows from the identity
\[
{1\over S_k-\lambda}f=\{l_k(\lambda,x)\}\ast f.\tag*{\qed}
\]
\renewcommand{\qed}{}
\end{proof}

\begin{corollary}
If $E_k(\lambda)\ne 0$, then
\begin{equation*}
{1\over(S_k-\lambda)^m}=\left\{{1\over(m-1)!}{\partial^{m-1}
\over\partial\lambda^{m-1}}l_k(\lambda,x)\right\}\ast.
\end{equation*}
\end{corollary}

\section{Heaviside algorithm for solving nonlocal Cauchy problems\\ for
Dunkl operators}\label{7}

Now we are to apply the elements of the operational calculus developed
in the previous section to ef\/fective solution of nonlocal Cauchy
boundary value problems of the form
\begin{equation}\label{Cauchy}
P(D_k)u=f,\qquad \Phi(D_k^ju)=\alpha_j,\qquad j=0,1,2,\dots,\deg P-1,
\end{equation}
with given $\alpha_j\in{\mathbb C}$.

To this end we extend the classical Heaviside algorithm, which is
intended for solving initial value problems for ordinary linear
dif\/ferential equations with constant coef\/f\/icients to the case of Dunkl
functional-dif\/ferential equations.

The extended Heaviside algorithm starts with the algebraization of
problem~(\ref{Cauchy}). It reduces the problem to a single algebraic
equation of the f\/irst degree in ${\cal M}_k$.

Let $P(\lambda)=a_0\lambda^m+a_1\lambda^{m-1}+\cdots+a_{m-1}\lambda
+a_m$ be a given polynomial of $m$-th degree, i.e.\ with $a_0\ne 0$.

The consecutive steps of the algorithm are the following:

1) Factorize $P(\lambda)$ in ${\mathbb C}$ to
\[
P(\lambda)=a_0(\lambda-\mu_1)^{\varkappa_1}(\lambda-\mu_2)
^{\varkappa_2}\cdots(\lambda-\mu_s)^{\varkappa_s},
\]
where $\mu_1,\mu_2,\dots,\mu_s$ are the distinct zeros of $P(\lambda)$
and $\varkappa_1,\varkappa_2,\dots,\varkappa_s$ are their
corresponding multiplicities.

2) Represent each of the terms of the equation by  the algebraic
Dunkl operator $S_k$. This is done by the formulae
\[
D_k^ju=S_k^ju-S_k^{j-1}\alpha_0-S_k^{j-2}\alpha_1-\cdots-S_k\alpha_{j-2}
-\alpha_{j-1}  ,\qquad j=1,2,\dots,m.
\]
Thus we obtain the following equation in ${\cal M}_k$:
\begin{equation*}
P(S_k)u=f+Q(S_k),\qquad \deg Q<\deg P,
\end{equation*}
with
\[
Q(S_k)=\sum_{j=0}^{m-1}\sum_{l=0}^{m-j-1}a_j\alpha_lS_k^{m-j-l-1}=
\sum_{\mu=0}^{m-1}\left(\sum_{\nu=0}^{m-\mu-1}
a_{\nu}\alpha_{m-\mu-\nu-1}\right)S_k^{\mu}.
\]

3) Verify if $P(S_k)$ is a non-divisor of zero in ${\cal M}_k$ by
checking if $E_k(\mu_j)\ne 0$ for all $j=1,2,\dots,s$.

4) If $P(S_k)$ is a non-divisor of zero, then write the solution $u$
in ${\cal M}_k$:
\[
u={1\over P(S_k)}f+{Q(S_k)\over P(S_k)}.
\]

5) Expand ${1\over P(S_k)}$ and
${Q(S_k)\over P(S_k)}$ into partial fractions:
\begin{gather*}
{1\over P(S_k)}=\sum_{j=1}^s\sum_{l=1}^{\varkappa_j}
{A_{j,l}\over(S_k-\mu_j)^l},
\qquad
{Q(S_k)\over P(S_k)}=\sum_{j=1}^s\sum_{l=1}^{\varkappa_j}
{B_{j,l}\over(S_k-\mu_j)^l}.
\end{gather*}

6) Interpret the partial fractions as convolution operators
\begin{gather*}
{1\over S_k-\mu_j}=\left\{l_k(\mu_j,x)\right\}\ast={1\over
E_k(\mu_j)}\left\{j_{k-{1\over 2}}(i\mu_j x)+{\mu_j x\over
2k+1}j_{k+{1\over 2}}(i\mu_j x)\right\}\ast.
\\
{1\over(S_k-\mu_j)^l}=\left\{{1\over(l-1)!}\left.{\partial^{l-1}
\over\partial\lambda^{l-1}}l_k(\lambda,x)\right|_{\lambda=\mu_j}
\right\}\ast\ ,\quad l=2,3,\dots.
\end{gather*}

7) Write the convolutional representation
\[
u(x)=(G\ast f)(x)+R(x),\qquad{\rm where}\qquad G={1\over P(S_k)},\qquad
R={Q(S_k)\over P(S_k)}.
\]

\begin{example}
Let $P(\lambda)$ has only simple zeros $\mu_1,\mu_2,\dots,\mu_m$. Then
\[
{1\over P(S_k)}=\sum_{j=1}^m{1\over P'(\mu_j)}\cdot{1\over
S_k-\mu_j}=\left\{\sum_{j=1}^m{1\over
P'(\mu_j)}l_k(\mu_j,x)\right\}\ast
\]
and
\[
{Q(S_k)\over P(S_k)}=
\sum_{j=1}^m{Q(\mu_j)\over P'(\mu_j)}\cdot{1\over S_k-\mu_j}
=\left\{\sum_{j=1}^m{Q(\mu_j)\over P'(\mu_j)}l_k(\mu_j,x)\right\}\ast.
\]
Then the solution $u$ takes the functional form
\[
u(x)=\sum_{j=1}^m{1\over P'(\mu_j)}l_k(\mu_j,x)\ast f(x)+
\sum_{j=1}^m{Q(\mu_j)\over P'(\mu_j)}l_k(\mu_j,x).
\]
\end{example}

The result of this section can be summarized in the following

\begin{theorem}\label{mainnlcp}
The nonlocal Cauchy problem $($Definition {\rm \ref{nlcp}}$)$ for a Dunkl
equation $P(D_k)u=f$ has a unique solution in $C^{(m)}({\mathbb R})$,
$m=\deg P$, iff none of the zeros of the polynomial $P(\lambda)$ is a~zero of the indicatrix $E_k(\lambda)$, i.e.\ when
\[
\{\lambda:P(\lambda)=0\}\cap\{\lambda:E_k(\lambda)=0\}=\varnothing.
\]
\end{theorem}

\begin{remark}
The term ``nonlocal'' should not be understood literally. The assertion
of Theorem~\ref{mainnlcp} is true also when $\Phi$ is a Dirac
functional, i.e.\ $\Phi\{f\}=f(a)$ for $a\in{\mathbb R}$. For us the
most interesting is the case $\Phi\{f\}=f(0)$. Then
$E_k(\lambda)\equiv 1$ and from the theorem it follows that the
initial value problem
\[
P(D_k)u=f,\qquad
u(0)=\alpha_0, \qquad (D_ku)(0)=\alpha_1, \qquad \dots, \qquad \big(D_k^{n-1}u\big)(0)=\alpha_{n-1},
\]
always has a unique solution. We will use this fact in the following
section.
\end{remark}

\section[Mean-periodic functions for $D_k$ determined by a linear
functional and mean periodic solutions of Dunkl equations]{Mean-periodic functions for $\boldsymbol{D_k}$ determined by a linear\\
functional and mean periodic solutions of Dunkl equations}\label{4}

The notion of mean-periodic function for the dif\/ferentiation operator
${d\over dt}$, determined by a linear functional $\Phi$ in
$C({\mathbb R})$, is introduced by J.~Delsarte~\cite{De}:

A function $f\in C({\mathbb R})$ is said to be {\it mean-periodic with
respect to the functional $\Phi$} if it satisf\/ies identically the
condition
\begin{equation}\label{mpf}
\Phi_{\tau}\{f(t+\tau)\}=0.
\end{equation}

In order to def\/ine mean-periodic functions for the Dunkl operator $D_k$
we need to recall the def\/inition of the Dunkl translation (shift)
operators, introduced by M.~R\"osler \cite{Ro} and later studied in M.A.~Mourou and K.~Trim\`eche~\cite{M-T}. They are a class of
operators $M:C({\mathbb R})\to C({\mathbb R})$ commuting with $D_k$ in
$C^1({\mathbb R})$.

\begin{definition}
Let $f\in C({\mathbb R})$ and $y\in{\mathbb R}$. Then
$(T_k^yf)(x)=u(x,y)\in C^1({\mathbb R}^2)$ is the solution of the
boundary value problem
\begin{equation*}
D_{k,x}u(x,y)=D_{k,y}u(x,y),\qquad u(x,0)=f(x).
\end{equation*}
$T_k^y$ is called the {\it translation operator for the Dunkl operator}
$D_k$.
\end{definition}
Such a solution exists for arbitrary $f\in C({\mathbb R})$ and it has
the following explicit form (see e.g.~\cite{Ro,S-K}):
\begin{gather*}
T_k^yf(x) =
{\Gamma\left(k+{1\over 2}\right)\over\Gamma(k)\Gamma
\left({1\over 2}\right)}
\left[\int_0^{\pi}f_{\rm e}\left(\sqrt{x^2+y^2-2|xy|\cos t}\right)
h^{\rm e}(x,y,t)\sin^{2k-1}t\,dt
\right.
\\
\left.\phantom{T_k^yf(x) = }{}  + \int_0^{\pi}f_{\rm o}\left(\sqrt{x^2+y^2-2|xy|\cos t}\right)
h^{\rm o}(x,y,t)\sin^{2k-1}t\,dt
\right].
\end{gather*}
As usually, the subscripts ``${\rm e}$'' and ``${\rm o}$'' denote correspondingly the
even and the odd part of a~function:
$f_{\rm e}(x)={f(x)+f(-x)\over 2}$, $f_{\rm o}(x)={f(x)-f(-x)\over2}$.
As for $h^{\rm e}(x,y,t)$ and $h^{\rm o}(x,y,t)$, they denote respectively
\begin{gather*}
h^{\rm e}(x,y,t) = 1-\sign(xy)\cos t,\\
h^{\rm o}(x,y,t) = \begin{cases}
\displaystyle
{(x+y)(1-\sign(xy)\cos t)\over\sqrt{x^2+y^2-2|xy|\cos t}}
&\text{for }(x,y)\ne(0,0),\\
0&$\text{otherwise.}$
\end{cases}
\end{gather*}

\begin{lemma}\label{proptrans}
The translation operators satisfy the following basic relations:
\begin{alignat}{3}
& (i) &&  T_k^yf(x)=T_k^xf(y), & \label{xyyx}\\
& (ii)  && T_k^yT_k^zf(x)=T_k^zT_k^yf(x), & \\
& (iii)\quad && D_{k,x}T_k^yf(x)=T_k^yD_{k,x}f(x). & \label{dttd}
\end{alignat}
\end{lemma}

Proofs can be found in various publications, in particular, in our
paper~\cite{D-H-S}.

A natural extension of the notion of mean-periodic function for the
Dunkl operator is proposed by Ben~Salem and Kallel \cite{S-K}. Instead of
(\ref{mpf}) they use the condition
\begin{equation}\label{mpfd}
\Phi_y\{T_k^yf(x)\}=0
\end{equation}
to def\/ine {\it mean-periodic function $f$ for $D_k$ with respect to the
functional $\Phi$}. Here $T_k^y$ is the generalized translation
operator just def\/ined.

{\it The space of mean-periodic functions for the Dunkl operator $D_k$
with respect to a given functional} $\Phi$ will be denoted by ${\cal
P}_{\Phi}$. We skip the subscript $k$ for sake of simplicity.

\begin{lemma}\label{lkmpf}
If $f\in {\cal P}_{\Phi}$, then $L_kf\in {\cal P}_{\Phi}$.
\end{lemma}
\begin{proof}
Denote $\varphi(x)=\Phi_t\{T_k^tL_kf(x)\}$ and use the commutation
relation (\ref{dttd}) from Lemma \ref{proptrans}
$D_{k,x}T_k^yf(x)=T_k^yD_{k,x}f(x)$ to obtain
\[
D_k\varphi(x)=\Phi_t\{D_kT_k^tL_kf(x)\}=\Phi_t\{T_k^tD_kL_kf(x)\}=
\Phi_t\{T_k^tf(x)\}=0.
\]
Hence $\varphi(x)=C={\rm const}$. But
$\varphi(0)=\Phi_t\{T_k^tL_kf(0)\}=\Phi_t\{T_k^0L_kf(t)\}=
\Phi_t\{L_kf(t)\}=0$.
Hence $C=0$.
\end{proof}

Further we will be interested in the solvability of Dunkl
dif\/ferential-dif\/ference equations
\begin{equation}\label{ddde}
P(D_k)u=f
\end{equation}
with a polynomial $P$ in the space of the mean-periodic functions
${\cal P}_{\Phi}$, def\/ined by (\ref{mpfd}). We intend also to propose
an algorithm for obtaining such solutions.

To this end we are to develop an operational calculus for $D_k$ in
$C({\mathbb R})$ and to extend the Heaviside algorithm for it. The
following result plays a basic role in the application of this
algorithm for solution of Dunkl equations in mean-periodic functions.

\begin{theorem}\label{ideal}
The class of mean-periodic functions ${\cal P}_{\Phi}$ is an ideal in
the convolutional algebra $(C({\mathbb R}),\ast)$, i.e.\ if $f\in {\cal
P}_{\Phi}$ and $g\in C({\mathbb R})$, then $f\ast g\in {\cal
P}_{\Phi}$.
\end{theorem}
\begin{proof} Assume that $f\in {\cal P}_{\Phi}$, i.e.
\[
\Phi_t\big\{T_k^tf(x)\big\}=0.
\]
From Lemma \ref{lkmpf} it follows that $L_k^{n+1}f\in {\cal P}_{\Phi}$
for
$n=0,1,2,\dots$, i.e.
\[
\Phi_t\big\{T_k^tL_k^{n+1}f(x)\big\}=0.
\]
Since $L_kf=\{1\}\ast f$, then $L_k^{n+1}f=A_{k,n}\ast f$, where the
Dunkl--Appell polynomial $A_{k,n}$ is of degree exactly $n$. We have
\[
\Phi_t\big\{T_k^t(A_{k,n}\ast f)(x)\big\}=0
\]
and then we can assert that
\[
\Phi_t\big\{T_k^t(P\ast f)(x)\big\}=0
\]
for any polynomial $P$.
By an approximation argument it follows that
\[
\Phi_t\big\{T_k^t(g\ast f)(x)\big\}=0
\]
for arbitrary $g\in C({\mathbb R})$, i.e.\ that $g\ast f\in {\cal
P}_{\Phi}$.
\end{proof}

\begin{corollary}
Let $M:C({\mathbb R})\to C({\mathbb R})$ be an arbitrary multiplier of
the algebra $(C({\mathbb R}),\ast)$. Then $M({\cal P}_{\Phi})\subset
{\cal P}_{\Phi}$, i.e.\ the restriction of $M$ to ${\cal P}_{\Phi}$ is
an inner operator in ${\cal P}_{\Phi}$.
\end{corollary}
\begin{proof}
Let $f\in {\cal P}_{\Phi}$. According to Theorem \ref{charact},
$Mf=D_k(m\ast f)$ with $m=M\{1\}$, then, by Theorem \ref{ideal},
$m\ast f\in {\cal P}_{\Phi}\cap C^1({\mathbb R})$. Then $D_k(m\ast
f)\in {\cal P}_{\Phi}$, i.e.\ $f\in {\cal P}_{\Phi}$ implies
$Mf\in {\cal P}_{\Phi}$.
\end{proof}

In the sequel we study the problem of solution of Dunkl equations in
mean-periodic functions determined by a linear functional.

\begin{theorem}\label{mps}
A function $u\in {\cal P}_{\Phi}\cap C^{(m)}({\mathbb R})$ is a
solution of the Dunkl equation $P(D_k)u=f$, with $f\in{\cal P}_{\Phi}$ iff $u$ is a
solution of the homogeneous nonlocal Cauchy problem
\[
P(D_k)u=f,\qquad \Phi\{D_k^ju\}=0,\qquad j=0,1,2,\dots,m-1,\qquad m=\deg P.
\]
\end{theorem}
\begin{proof}
The condition $f\in {\cal P}_{\Phi}$ is necessary for the existence of
a solution $u\in {\cal P}_{\Phi}$. Assume that a~function $u\in {\cal
P}_{\Phi}\cap C^{(m)}({\mathbb R})$ is a solution of the Dunkl
equation $P(D_k)u=f$.
Then mean-periodic are all the functions $D_k^ju$,
$j=0,1,2,\dots,m-1$, i.e.
\begin{equation}\label{mp1}
\Phi_y\big\{T_k^yD_k^ju(x)\big\}=0,
\end{equation}
since the operator $Af(x)=\Phi_y\{T_k^yf(x)\}$ commutes with $D_k$
(Dimovski, Hristov and Sif\/i~\cite{D-H-S}).
For $x=0$ from~(\ref{mp1}) we get
\[
\Phi_y\big\{T_k^yD_k^ju(0)\big\}=0.
\]
But $T_k^yD_k^ju(0)=T_k^0D_k^ju(y)$ ((\ref{xyyx}), Lemma~\ref{proptrans}) and hence
\begin{equation}\label{mp2}
\Phi\big\{D_k^ju\big\}=0,\qquad j=0,1,2,\dots,m-1.
\end{equation}
In order to prove that a solution $u$ of $P(D_k)u=f$ with $f\in
{\cal P}_{\Phi}$, which satisf\/ies conditions (\ref{mp2}), is a
mean-periodic function, we consider the function
\[
v=\Phi_y\big\{T_k^yu(x)\big\}=Au.
\]
Since the operator $A$ commutes with $D_k$, then applying it on the
equation $P(D_k)u=f$, we get $P(D_k)v=0$ due to $Af=0$. It remains to
f\/ind the initial values $D_k^jv(0)$, $j=0,1,2,\dots,m-1$:
\[
D_k^jv(0)=AD_k^ju(0)=\Phi_y\{T_k^yD_k^ju(0)\}=\Phi_y\{T_k^0D_k^ju(y)\}
=\Phi_y\{D_k^ju(y)\}=0.
\]
At the end of the previous section we have seen that the initial value
problem $P(D_k)v=0$, $D_k^jv(0)=0$, $j=0,1,2,\dots,m-1$, has only the
trivial solution $v(x)=0$. Thus we proved that $\Phi_y\{T_k^yu\}=0$,
i.e.~$u$ is mean-periodic.
\end{proof}

Now we can use operational calculus method for solving nonlocal Cauchy
problems for Dunkl equations to f\/ind explicitly the mean-periodic
solutions of such equations.

To this end, we are to solve the homogeneous nonlocal Cauchy boundary
value problem
\begin{equation}\label{homog}
P(D_k)u=f,\qquad\Phi\big\{D_k^ju\big\}=0,\qquad j=0,1,2,\dots,m-1,
\end{equation}
with $f\in {\cal P}_{\Phi}$.

In the ring ${\cal M}_k$ of the multiplier fractions it reduces to
the single algebraic equation for $u$
\begin{equation}\label{mp3}
P(S_k)u=f.
\end{equation}
As we have seen in Section \ref{6}, $P(S_k)$ is a non-divisor of zero
in ${\cal M}_k$ if\/f none of the zeros of the polynomial $P(\lambda)$
is a zero of the Dunkl indicatrix $E_k(\lambda)$. If $P(S_k)$ is a~divisor of zero, then, in order to ensure the existence of solution of~(\ref{mp3}) and thus of~(\ref{homog}), additional restrictions on $f$
should be imposed. This is the so called {\it resonance case}, which
we will not treat here.

Thus, let $P(S_k)$ be a non-divisor of zero in ${\cal M}_k$, i.e.\
$\{\lambda:P(\lambda)=0\}\cap\{\lambda:E_k(\lambda)=0\}=\varnothing$.
Then the formal solution of (\ref{mp3}) in ${\cal M}_k$
\[
u={1\over P(S_k)}f
\]
can be written in explicit functional form. Using the extended
Heaviside algorithm of Section~\ref{7}, we represent
$ {1\over P(S_k)}$ as a convolutional operator
\[
{1\over P(S_k)}=\{G(x)\}\ast.
\]
Then
\begin{equation*}
u=G\ast f
\end{equation*}
is the desired mean-periodic solution of the Dunkl equation
$P(D_k)u=f$. The verif\/ication is straightforward. Indeed, $G\ast f\in
{\cal P}_{\Phi}$ according to Theorem~\ref{ideal}, since $f\in {\cal
P}_{\Phi}$.

Our considerations of the problem for solving Dunkl equations in
mean-periodic functions can be summarized in the following

\begin{theorem}\label{mpsol}
A Dunkl equation $P(D_k)u=f$ with $f\in {\cal P}_{\Phi}$ has a unique
solution in ${\cal P}_{\Phi}$ iff none of the zeros of the polynomial
$P(\lambda)$ is a zero of the Dunkl indicatrix
\[
E_k(\lambda)=\Phi\left\{j_{k-{1\over 2}}(i\lambda x)+{\lambda x\over
2k+1}j_{k+{1\over 2}}(i\lambda x)\right\}.
\]
\end{theorem}

In the end, it is possible the Duhamel principle to be extended to the
problem for solving Dunkl equations in mean-periodic functions.

\begin{theorem}\label{hommps}
Let $H(x)$ be the solution of the homogeneous nonlocal Cauchy problem
$P(D_k)H{=}1$, $\Phi\{D_k^jH\}=0$, $j=0,1,2,\dots,m-1$. Then
\[
u=D_k(H\ast f)
\]
is a mean-periodic solution of the Dunkl equation $P(D_k)u=f$ with
$f\in {\cal P}_{\Phi}$.
\end{theorem}

\subsection*{Acknowledgments}
 The authors are very grateful to the editors
and to the referees for the constructive and valuable comments and
recommendations.

\pdfbookmark[1]{References}{ref}
\LastPageEnding

\end{document}